\documentclass[12pt, twoside]{amsart}
\usepackage[english]{babel}
\theoremstyle{plain}
\addtocounter{page}{1}
\usepackage{graphicx, color}
\usepackage{cancel}

\newtheorem{theorem}{Theorem}[section]
\newtheorem{remark}{Remark}[section]

\newtheorem{corollary}{Corollary}[section]

\title[Radial $p$-Laplacian supercritical Neumann problems]{Radial positive solutions for
$p$-Laplacian supercritical Neumann
problems\smallskip
\\
Soluzioni radiali positive di problemi di Neumann supercritici governati dal $p$-laplaciano
}
\author{F. Colasuonno}
\address{Francesca Colasuonno\newline\indent
Dipartimento di Matematica, 
Alma Mater Studiorum Universit\`a di Bologna
\newline\indent
piazza di Porta S. Donato 5, 40126 Bologna, Italia}
\email{francesca.colasuonno@unibo.it}

\author{B. Noris}
\address{Benedetta Noris
\newline \indent Laboratoire Ami\'enois de Math\'ematique Fondamentale et Appliqu\'ee\newline\indent
Universit\'e de Picardie Jules Verne\newline\indent
33 rue Saint- Leu, 80039 AMIENS, France}
\email{benedetta.noris@u-picardie.fr}

\begin{document}

\begin{abstract}
{This paper deals with existence and multiplicity of positive solutions for a quasilinear problem with Neumann boundary conditions, set in a ball. The problem admits at least one constant non-zero solution and it involves a nonlinearity that can be supercritical in the sense of Sobolev embeddings. The main tools used are variational techniques and the shooting method for ODE's. These results are contained in \cite{BF, ABF}.
}

\medskip \noindent {\sc{Sunto.}} 
{In questo lavoro trattiamo l'esistenza e la molteplicit\`a di soluzioni positive per un probelma quasili\-ne\-are ambientato in una palla, con condizioni al bordo di Neumann. Il problema ammette almeno una soluzione costante non nulla e coinvolge una nonlinearit\`a che pu\`o essere supercritica nel senso delle immersioni di Sobolev. I principali strumenti usati nello studio di tale problema sono tecniche variazionali e il metodo di shooting per le EDO. Questi risultati sono contenuti in \cite{BF, ABF}.
}


\end{abstract}
\maketitle
\section{Introduction}
For $1<p<\infty$, we consider the following quasilinear Neumann problem
\begin{equation}\label{Pg}
\begin{cases}
-\Delta_p u+u^{p-1}=g(u)\quad&\mbox{in }B_R,\\
u>0\quad&\mbox{in }B_R,\\
\partial_\nu u=0\quad&\mbox{on }\partial B_R,
\end{cases}
\end{equation}
where $\Delta_p u:=\mathrm{div}(|\nabla u|^{p-2}\nabla u)$ denotes the $p$-Laplace operator, $B_R\subset \mathbb R^N$ is the ball of radius $R$ centered at the origin, $N\ge1$, and $\nu$ is the outer unit normal of $\partial B_R$.
In \cite{BF,ABF}, we investigate the existence of non-constant solutions of \eqref{Pg} under very mild assumptions on the nonlinearity $g$, allowing in particular for supercritical growth in the sense of Sobolev embeddings. We observe that, differently from Dirichlet supercritical problems, in the case of Neumann boundary conditions there is not a Pohozaev-type obstruction to the existence of non-zero solutions, so the natural question that arises is whether the problem admits any {\it non-constant} solutions. 

We will show that the situation changes drastically depending on $p>1$. We start with considering the case $p\ge 2$. 

Let us first introduce the assumptions on the nonlinearity. In \cite{BF}, we assume that $g:[0,\infty)\to\mathbb R$ is of class $C^1([0,\infty))$ and satisfies the following hypotheses\footnote{In \cite{BF}, the hypothesis $(g_0)$ requires the limit in 0 to belong to $[0,1)$ instead of $(-\infty,1)$. Nevertheless, that assumption can be weakened as stated here, because it is always possible to modify $g(s)$ into $\tilde g(s):=g(s)+ms^{p-1}$ for a suitable $m>0$ such that $\lim_{s\to 0^+}\tilde g(s)/s^{p-1}\in[0,1)$, and study the equivalent problem $-\Delta_p u+(m+1)u^{p-1}=\tilde g(u)$. We observe in passing that the constant $m$ can be also adjusted in such a way to deal with a non-negative and non-decreasing $\tilde g$.}

\begin{minipage}[l]{.55\textwidth}
\begin{itemize}
\item[$(g_0)$] $\lim_{s\to 0^+}\frac{g(s)}{s^{p-1}}\in (-\infty,1)$;
\item[$(g_\infty)$] $\liminf_{s\to\infty}\frac{g(s)}{s^{p-1}}\in(1,\infty]$;
\item[$(g_{u_0})$] $\exists$ $u_0>0$ such that $g(u_0)=u_0^{p-1}$ and 
$$
g'(u_0)>
\begin{cases}
(p-1)u_0^{p-2}\quad&\mbox{if }p>2,\\
\lambda_{2}^{\textnormal{rad}}+1&\mbox{if }p=2,
\end{cases}$$
where $\lambda_2^{\mathrm{rad}}$ denotes the second radial eigenvalue of Neumann Laplacian $-\Delta$.
\end{itemize}
\end{minipage}\quad
\begin{minipage}[r]{.45\textwidth}
\includegraphics[width=6cm]{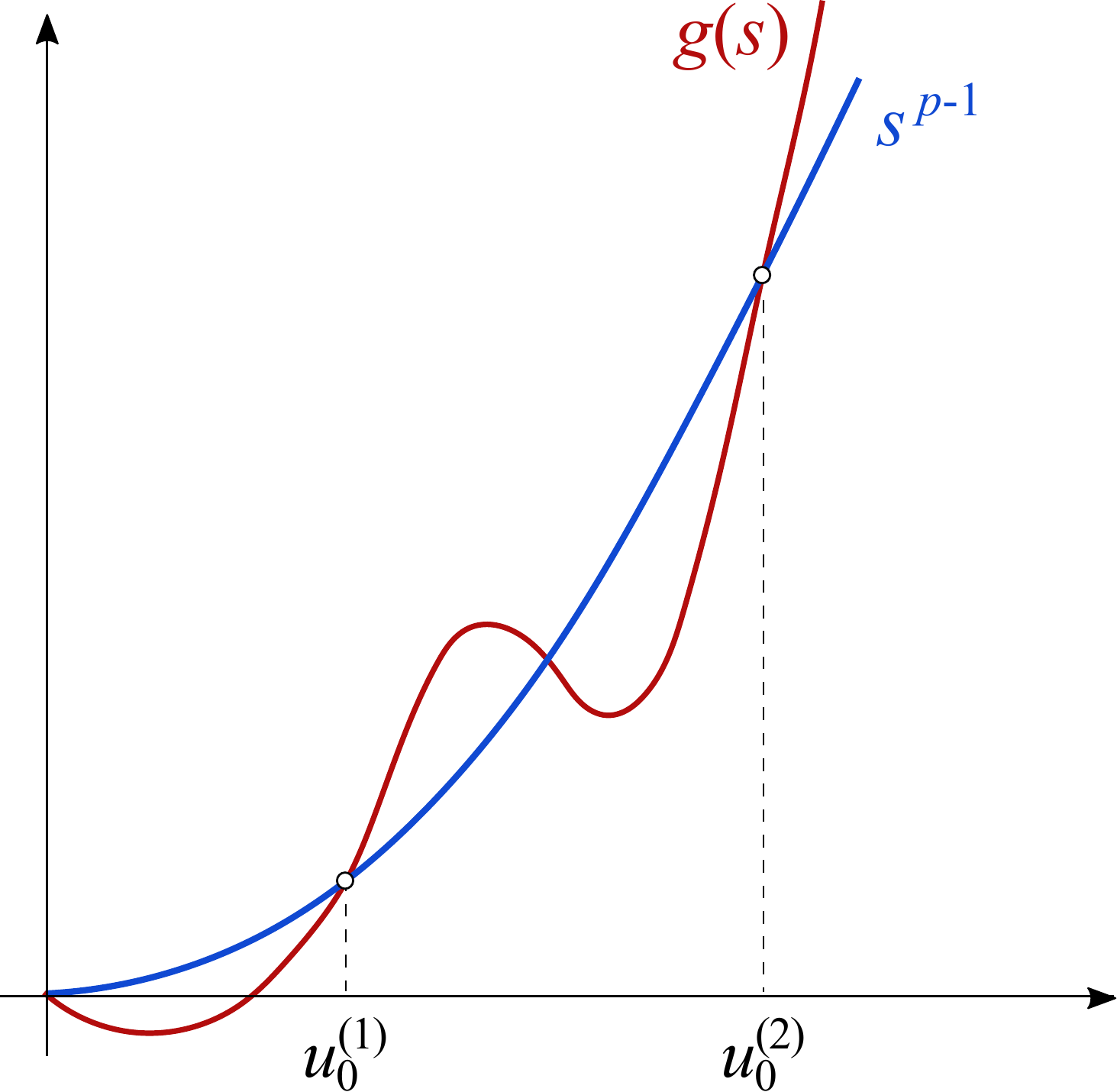}
\end{minipage}
\smallskip 

The prototype nonlinearity $g$ is the pure power $s^{q-1}$ for $q>p$. 
\smallskip 

As an immediate consequence of $(g_0)$, $g(0)=0$. Moreover, by $(g_\infty)$, the nonlinearity $g$ can be taken Sobolev-supercritical. We further remark that, by the regularity of $g$ and by $(g_0)$ and $(g_\infty)$, we immediately have the existence of an intersection point $u_0>0$ between $g$ and the power function $s^{p-1}$, with $g'(u_0)\ge (s^{p-1})'(u_0)= (p-1)u_0^{p-2}$. Hence, when $p>2$, condition $(g_{u_0})$ is only needed to prevent the situation in which $g$ is tangent to $s^{p-1}$ at $u_0$. While for $p=2$, the condition required at $u_0$ is stronger, being $\lambda_2^{\mathrm{rad}}>0$. In both cases, $p>2$ and $p=2$, conditions $(g_0)$ and $(g_\infty)$ are enough to prove the existence of a radial solution to \eqref{Pg} of minimax-type, while $(g_{u_0})$ is needed to prove that the solution found is non-constant.
We finally observe that, due to the existence of $u_0>0$ for which $g(u_0)=u_0^{p-1}$, problem \eqref{Pg} admits at least the constant solution $u\equiv u_0$. 

The main result in \cite{BF} reads as follows. 

\begin{theorem}[{\cite[Theorem 1.3]{BNW} for $p=2$,\cite[Theorem 1.1]{BF} for $p>2$}]\label{thm:intro}
Let $p\geq2$ and let $g$ satisfy the hypotheses above.
There exists a non-constant, radial, radially non-decreasing solution of \eqref{Pg}.
In addition, if $u_{0,1},\ldots,u_{0,n}$ are $n$ different positive constants satisfying $(g_{u_0})$, then \eqref{Pg} admits $n$ different non-constant, radial, radially non-decreasing solutions.
\end{theorem}

Let us now spend a few words on the techniques used to prove Theorem \ref{thm:intro}.
Since the equation in \eqref{Pg} is possibly supercritical, the energy associated to the problem might not be well defined in the whole of $W^{1,p}(B_R)$, and so, a priori, variational methods cannot be used to study this problem. Nevertheless, we take advantage of the idea proposed by Serra and Tilli in \cite{ST} and work in the cone of non-negative, radial, radially non-decreasing functions  
\begin{equation}\label{cone}
\mathcal C:=\Big\{u\in W^{1,p}_{\mathrm{rad}}(B_R)\,:\, u\ge0,\,u(r)\le u(s) \mbox{ for all }0<r\le s\le R\Big\},
\end{equation}
where with abuse of notation we write $u(|x|):=u(x)$. The main reason for working in this set is that {\it all solutions of \eqref{Pg} belonging to $\mathcal C$ are a priori bounded in $W^{1,p}(B_R)$ and in $L^\infty(B_R)$.}
By the way, the cone $\mathcal C$ {\it has empty interior in the $W^{1,p}$-topology}, so in general, if we define the associated energy functional $I_\mathcal C:\mathcal C\to\mathbb R$, a function $u$ such that 
$$I'_\mathcal C(u)[\varphi]=0\quad\mbox{for all }\varphi\in\mathcal C$$
is not a weak solution of \eqref{Pg}.
The strategy used in \cite{BNW,BF} to overcome this difficulty is based on the truncation method. A sketch of the proof of Theorem \ref{thm:intro} is given in Section~\ref{pf-BF}, see also \cite{C}.

\bigskip

In \cite{ABF}, we consider problem \eqref{Pg} for every $p>1$. We require slightly different conditions to $g$. Namely, we assume less regularity, $g\in C([0,\infty)\cap C^1((0,\infty))$, and suppose that it satisfies the following assumptions 
\begin{itemize}
\item[$(g_0)'$] $\lim_{s\to 0^+}\frac{g(s)}{s^{p-1}}\in(-\infty,1]$;
\item[$(g_{\mathrm{eq}})$] $g(s)-s^{p-1}\,\begin{cases}
<0\quad&\mbox{if }\;0<s<1\\
=0&\mbox{if }\;s=1\\
>0&\mbox{if  }\;s>1;\\
\end{cases}
$
\item[$(g_1)$] there exists $C_1\in[0,\infty]$ such that $\lim_{s\to1}\frac{g(s)-s^{p-1}}{|s-1|^{p-2}(s-1)}=C_1$.
\end{itemize}
We note that $(g_{\mathrm{eq}})$ means that $g$ intersects only {\it once} the power $s^{p-1}$ at a point $u_0$, which without loss of generality is taken equal to 1. We believe that this condition can be weakened in order to allow more than one intersection between $g$ and $s^{p-1}$. 
Furthermore, we observe that while the assumption in zero (i.e., $(g_0)'$) is just slightly more general than before (i.e., $(g_0)$), we have replaced $(g_{u_0})$ with $(g_1)$. Condition $(g_1)$ is implied by the regularity of $g$ when $1<p\le 2$. Indeed, since $g$ is of class $C^1$ at 1, hypothesis $(g_1)$ holds automatically with $C_1\in[0,\infty)$ for $p=2$, and with $C_1=0$ for $p<2$.
The only case in which the existence of the limit in $(g_1)$ is not implied by the regularity of $g$ (and consequently $(g_1)$ is really an additional assumption) is when $p>2$ and $g'(1)=0$. Furthermore, we observe that for $p>2$, condition $(g_{u_0})$ required in \cite{BF} is stronger than $(g_1)$, since  $(g_{u_0})$ (for $u_0=1$) implies $(g_1)$ with $C_1=\infty$.
 
With this set of hypotheses, the prototype nonlinearity can be taken also of the form 
$$g(s)=s^{q-1}+s^{p-1}-s^{r-1}\quad\mbox{with }p\le r<q,$$
so that in general the prototype equation becomes 
$$-\Delta_p u+u^{r-1}=u^{q-1}\quad\mbox{in }B_R.$$ 
We further remark that in \cite{ABF} it is also treated the case set in an {\it annular domain}. Since the arguments are similar to the ones for the ball, for the sake of simplicity we present here only the case of the ball. The main result in \cite{ABF} is the following.

\begin{theorem}[Theorems 1.2 and 1.4 of \cite{ABF}]\label{thm:ABF}
Let $p>1$ and $\lambda_{k}^{\mathrm{rad}}$ denote the $k$-th radial eigenvalue of $-\Delta_p$ with Neumann boundary conditions for any integer $k\ge1$. If $g$ satisfies $(g_0)'$-$(g_1)$, then the following implications hold.  
\begin{itemize}
\item[(i)] If $C_1>\lambda_{k+1}^{\mathrm{rad}}$, then \eqref{Pg} admits at least $k$ different non-constant radial solutions $u_1,\dots,u_k$. Furthermore, $u_j-1$ has exactly $j$ zeros for any $j=1,\dots,k$.
\item[(ii)] If $C_1=\infty$, then 
problem \eqref{Pg} admits infinitely many non-constant radial solutions.
\item[(iii)] If $C_1=0$, then\footnote{When the domain is an annulus $A(R_1,R_2)$, part (iii) of Theorem \ref{thm:ABF} reads as 

{\it If $C_1=0$, then for every integer $k\ge1$ and any $\varepsilon>0$ there exists $R_*(k,\varepsilon)>0$ such that  if $R_1<\varepsilon R_2$ and $R_2>R_*(k,\varepsilon)$, \eqref{Pg} admits at least $2k$ different non-constant radial solutions $u_1^\pm,\dots,u_k^\pm$.}

The oscillating behavior is the same as for the solutions in the ball. } 
for every integer $k\ge1$ there exists $R_*(k)>0$ such that  if $R>R_*(k)$, \eqref{Pg} admits at least $2k$ different non-constant radial solutions $u_1^\pm,\dots,u_k^\pm$. Furthermore,  $u_j^\pm-1$ has exactly $j$ zeros for any $j=1,\dots,k$.
\end{itemize}
\end{theorem}
Clearly, part (ii) of the previous theorem can be seen as an immediate consequence of part (i), being $C_1=\infty$ greater than every eigenvalue $\lambda_k^{\mathrm{rad}}$. 
Now, when $g(s)= s^{q-1}$  with $q > p$, the constant $C_1$ in condition $(g_1)$ specializes in 
$$
C_1=\begin{cases}
+\infty\quad&\mbox{if }p>2,\\
q-2\quad&\mbox{if }p=2,\\
0\quad&\mbox{if }1<p<2,
\end{cases}
$$
and consequently Theorem \ref{thm:ABF} becomes
\begin{corollary}[Corollary 1.5 of \cite{ABF}]\label{cor:proto}
Let $g(s)=s^{q-1}$ with $q>p$. 
\begin{itemize}
\item[(i)] If $p=2$ and $q-2>\lambda_{k+1}^{\mathrm{rad}}$ for some $k\ge1$, then \eqref{Pg} admits at least $k$ different non-constant radial solutions $u_1,\dots,u_k$. Furthermore, $u_j-1$ has exactly $j$ zeros for any $j=1,\dots,k$.
\item[(ii)] If $p>2$, then \eqref{Pg} admits infinitely many non-constant radial solutions. 
\item[(iii)] If $1<p<2$, then for every integer $k\ge1$ there exist $R_*(k)>0$ such that if $R>R_*(k)$, problem \eqref{Pg} admits at least $2k$ different non-constant radial solutions $u_1^\pm,\dots,u_k^\pm$. Furthermore, $u_j^\pm-1$ has exactly $j$ zeros for any $j=1,\dots,k$.
\end{itemize}
\end{corollary}

Part (i) of the previous corollary is the same result as in \cite{bonheure2016multiple} (see Theorem~\ref{thm:BGT}-(i) below), but the proof techniques are completely different. We also observe that the condition on the exponent, i.e., $q>2+\lambda_{k+1}^{\mathrm{rad}}$, can be also read in terms of the radius $R$ of the ball: since the eigenvalues $\lambda_{k}^{\mathrm{rad}}=\lambda_{k}^{\mathrm{rad}}(R)$ are decreasing in $R$, keeping $q$ fixed, we can increase the radius $R$ in order to have the condition satisfied. In this way, the assumption in (i) becomes much more akin to the one in (iii).  Moreover, from (iii) we can see that for $1<p<2$  a completely different behavior appears: non-constant solutions with the same oscillatory behavior come in couples as soon as the radius of the domain overcomes a certain threshold. 

This note is organized as follows. In Section~\ref{pf-BF}, we sketch the proof of Theorem \ref{thm:intro}, while in Section~\ref{simulations} we collect some comments, pre-existing results, and numerical simulations to get further insights into the features of the solutions when $p=2$ and $p>2$. In Section~\ref{Sec4} we deal with the proof of Theorem~\ref{thm:ABF}-(i), and we conclude the paper by illustrating, in Section~\ref{Sec5}, the main reasons why the result for $p<2$ differs so much from the ones for $p\ge2$, through the guidelines of the proof of Theorem \ref{thm:ABF}-(iii) and the description of some numerical results.  
 
\section{Proof of Theorem \ref{thm:intro}}\label{pf-BF}
We sketch here the proof of Theorem \ref{thm:intro}. As mentioned in the Introduction, we restrict ourselves to the cone $\mathcal C$ of non-negative, radial, radially non-decreasing $W^{1,p}$-functions defined in \eqref{cone}, where it is possible to find a priori estimates on the solutions of \eqref{Pg}. We split the proof of the theorem into four steps. 

\underline{\it Step 1.} {\it(Truncation)} Thanks to the a priori estimates, we can {\it truncate the nonlinearity} $g$ and redefine it at infinity, in order to deal with a {\it subcritical nonlinearity}. In this way, we end up with a new {\it truncated problem} with the property that all solutions of the truncated problem belonging to $\mathcal C$ solve also the original problem \eqref{Pg}. 

\underline{\it Step 2.} {\it(Existence)} The energy functional $I$ associated to the truncated problem is well defined in the whole of $W^{1,p}(B_R)$, hence we can now apply variational methods. We need to find a critical point of $I$ which belongs to $\mathcal C$. To this aim, we prove that a mountain pass-type theorem holds for $I$ inside the cone $\mathcal C$. The main difficulty here is the construction of a descending flow that preserves $\mathcal C$, cf. \cite[Lemmas 3.7-3.8]{BF}. When $p>2$, this step presents the additional technical difficulty of proving the existence of a local Lipschitz vector field that preserves the cone $\mathcal C$, see \cite[Lemmas 3.4-3.6]{BF}.

\underline{\it Step 3.} {\it(Non-constancy)} We want to prove that the solution found is nonconstant. To this aim, we further restrict our cone, working in a subset of $\mathcal C$ in which the only constant solution of \eqref{Pg} is the positive constant $u_0$ defined in $(g_{u_0})$. In this set, we build an admissible curve along which the energy is lower than the energy of the constant $u_0$, which gives immediately that the minimax solution found (whose energy is such that $I(u)=\min_{\gamma\in\Gamma}\max_{t\in[0,1]}I(\gamma(t))$, where $\Gamma$ is the set of admissible curves) is not identically equal to $u_0$. More precisely, let $\phi_2$ be the second eigenfunction of the Neumann $p$-Laplacian. Via second-order Taylor expansion of $I$, we prove that for every $s\in(-\varepsilon,\varepsilon)\setminus\{0\}$
 $$
\begin{aligned}
I(t(s)(u_0+s\phi_2)) &- I(u_0)=\\
&\hspace{-1.3cm}\begin{cases}
\dfrac{s^2}{2}\displaystyle{\int_{B_R}}\Big\{\cancel{|\nabla u_0|^{p-2}|\nabla \phi_2|^2}+[(p-1)u_0^{p-2}-g'(u_0)]\phi_2^2\Big\}\,dx+o(s^2)<0\;&(p>2),\smallskip \\
\dfrac{s^2}{2}\displaystyle{\int_{B_R}}\Big\{|\nabla \phi_2|^2+ [1-g'(u_0)]\phi_2^2\Big\}dx+o(s^2)<0&(p=2),
\end{cases}
\end{aligned}
$$
where $t(s)$ is a suitable continuous function. We stress that the inequality signs in the above computation, both for $p>2$ and for $p=2$, are due to condition $(g_{u_0})$. This makes apparent the reason why we need to require different conditions for $p>2$ and $p=2$. Now, to get the admissible curve $\gamma\in\Gamma$ along which the energy is lower than $I(u_0)$, it is enough to rescale suitably the curve $s\mapsto t(s)(u_0+s\phi_2)$.
Finally, we observe here that this part of the proof uses heavily the $C^2$-regularity of the energy functional $I$, thus it cannot be generalized to the case $1<p<2$.

\underline{\it Step 4.} {\it(Multiplicity)} If there is more than one constant $u_0$ satisfying condition $(g_{u_0})$, we take advantage of the fact that, since we work in the restricted cone containing exactly one constant solution, we automatically localize each minimax solution. This allows us to prove the multiplicity result stated in Theorem \ref{thm:intro}, by simply repeating the same argument in each cone restricted about each $u_{0,i}$.
  
\section{Some comments on the case $p\ge2$}\label{simulations}

\begin{figure}[h!t]
\includegraphics[width=0.9\textwidth]{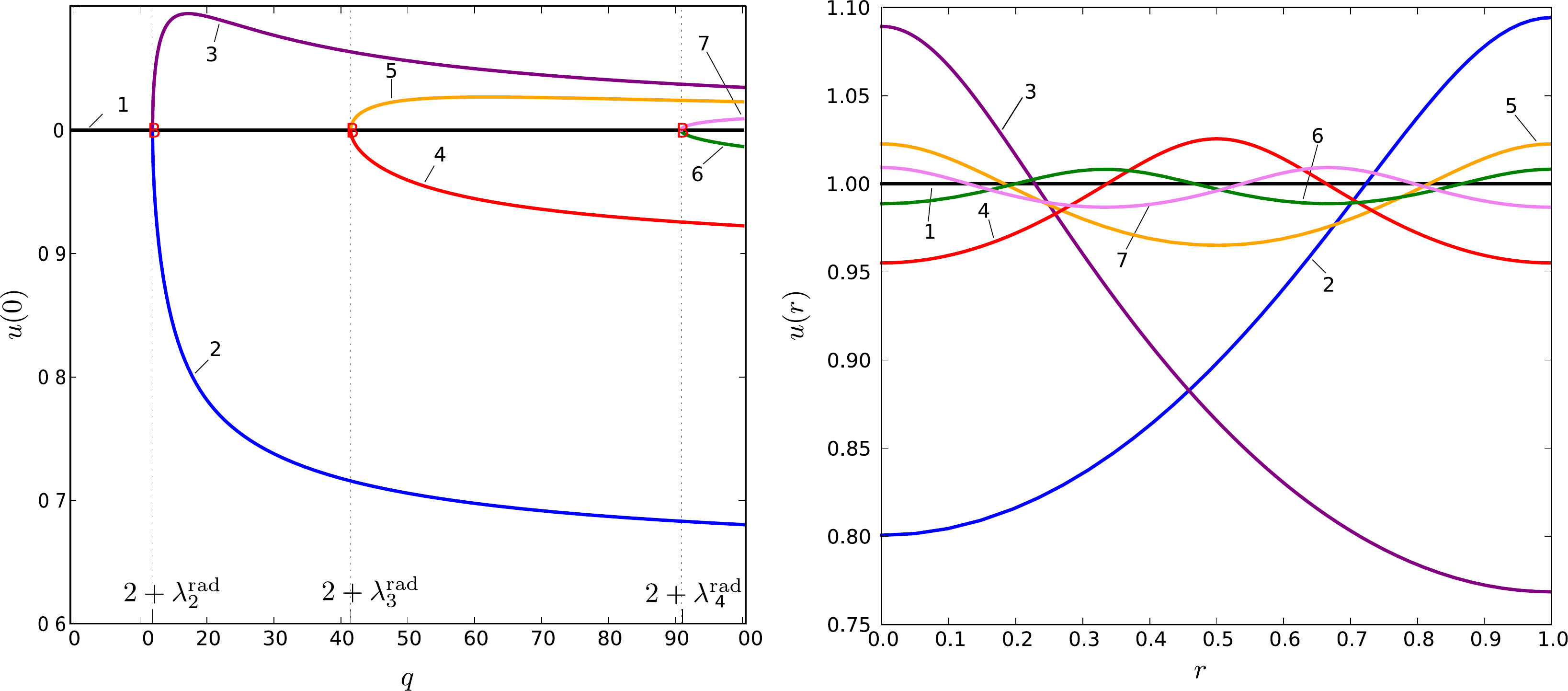}
\caption{The first three bifurcation branches for problem \eqref{Pg} in the case $N=1$, $R=1$, $p=2$, $g(s)=s^{q-1}$. {\it On the left}: bifurcation diagram $u(0)$ as function of $q$. {\it On the right}: solutions belonging to the first three branches. The color of each solution in the right plot corresponds to the color of the branch it belongs to in the left plot. More precisely, the numbers along the branches in the left plot are located in correspondence with the initial condition $u(0)$ of the solution represented in the right plot.}
\label{fig:p2}
\end{figure}

\begin{figure}[h!t]
\begin{center}
\includegraphics[width=\textwidth]{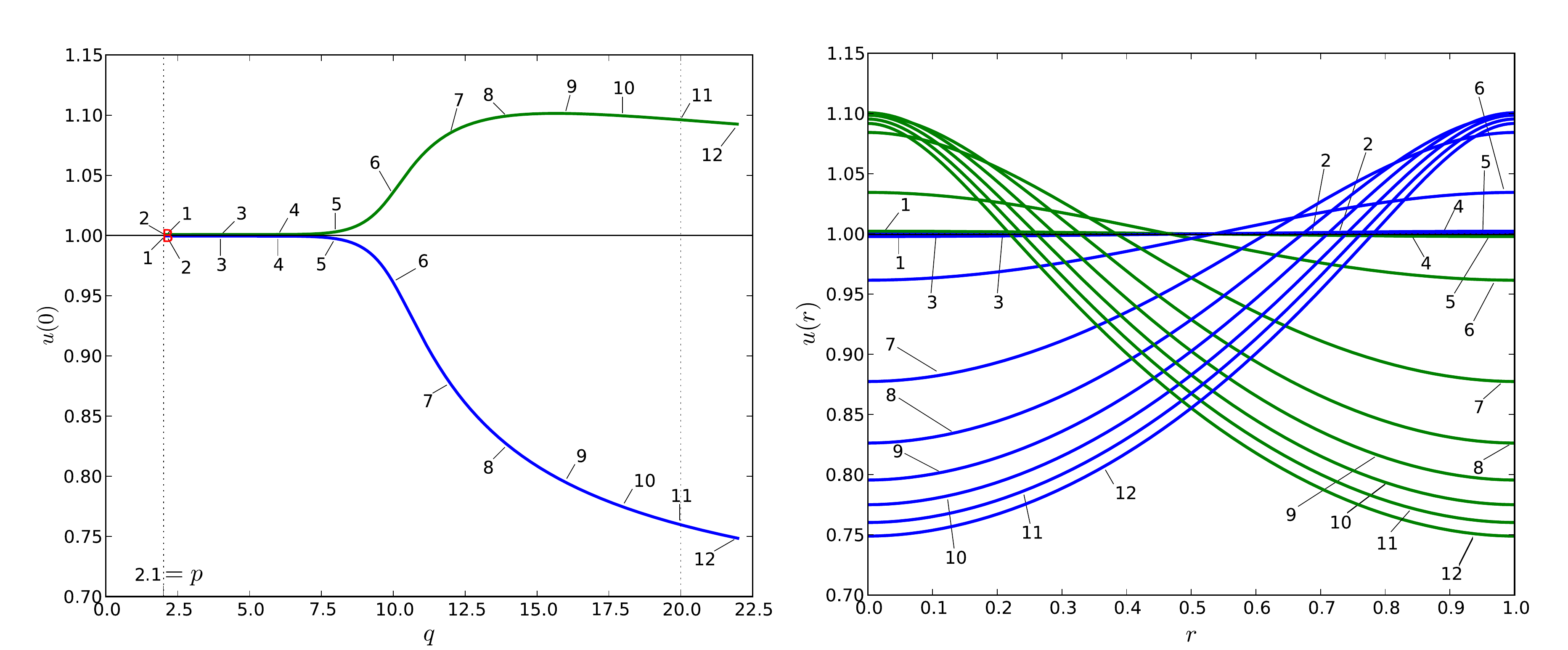}
\caption{Partial bifurcation diagram for problem \eqref{Pg} in the case $N=1$, $R=1$, $p=2.1$, and $g(s)=s^{q-1}$. {\it On the left}: bifurcation diagram $u(0)$ as function of $q$. The first two branches of solutions bifurcating at $q=p>2$; the green one is the branch of decreasing solutions, the blue one is the branch of increasing solutions. {\it On the right}: solutions belonging to the first branches. Blue (increasing) solutions belong to the blue branch in the left plot, green (decreasing) solutions belong to the green branch.}
\label{fig:ps2}
\end{center}
\end{figure}

From Step 3. above, one could get the impression that condition $(g_{u_0})$ is only a technical {\it ad hoc} assumption imposed on $g$ in order to let the machinery of the proof work fine. Actually, with reference to the bifurcation diagrams in Figures \ref{fig:p2} and \ref{fig:ps2}, one can see that the values $q=p$ for $p>2$ and $q=2+\lambda_2^{\mathrm{rad}}$ for $p=2$, involved in condition $(g_{u_0})$ when $g(s)=s^{q-1}$, arise naturally from the problem. Despite this, one should be aware that it has been proved in \cite{bonheure2016multiple} that, for $p=2$ and $N\ge 3$, the value $2+\lambda_2^{\mathrm{rad}}$ is not sharp.  
\medskip

Let us first comment the case $p=2$. We notice that, in the semilinear case, condition $(g_{u_0})$ involves the second radial eigenvalue of $-\Delta$ with Neumann boundary conditions. This is coherent with the result in \cite{bonheure2016multiple}, where the authors show that a bifurcation phenomenon underlies the existence result, at least in the case of the prototype nonlinearity $g(s)=s^{q-1}$. They prove that at $q=2+\lambda_{k+1}^{\textnormal{rad}}$, $k\geq1$, a new branch of solutions bifurcates from the constant branch $u\equiv u_0=1$. 
\begin{theorem}[\cite{bonheure2016multiple}]\label{thm:BGT}
Let $p=2$, $g(s)=s^{q-1}$ with $q>2$, and $\lambda_{k}^{\mathrm{rad}}$ denote the $k$-th eigenvalue for the Neumann Laplacian for any integer $k\ge1$.   
\begin{itemize}
\item[(i)] If $q>2+\lambda_{k+1}^{\mathrm{rad}}$, there exist at least $k$ non-constant radial solutions $u_1,\dots,u_k$ of \eqref{Pg}. Furthermore, $u_j-1$ has exactly $j$ zeros for any $j=1,\dots,k$.
\item[(ii)] If $2^*>q>2+\lambda_{k+1}^{\mathrm{rad}}$ (where $2^*$ is the Sobolev critical exponent), there exist at least $2k$ non-constant radial solutions $u_1^\pm,\dots,u_k^\pm$ of \eqref{Pg}. Furthermore, $u_j^\pm-1$ has exactly $j$ zeros for any $j=1,\dots,k$.
\end{itemize}
\end{theorem}
This theorem was proved by means of the Crandall-Rabinowitz bifurcation technique in the parameter $q$. As already mentioned in the Introduction, part (i) of the previous theorem was also recovered in \cite[Corollary 1.5-(ii)]{ABF} via shooting method.

We present now some numerical simulations performed with the software AUTO-07p for problem \eqref{Pg} in dimension $N=1$, with $R=1$ and $g(s)=s^{q-1}$. 

In Figure~\ref{fig:p2}, we represent the first three bifurcation branches for this problem with $p=2$. The black line represents the constant solution $u\equiv1$; the branches bifurcate at points $q=2+\lambda_k^{\text{rad}}$, $k=2,3,4$. The solutions belonging to the lower part of the first branch are monotone increasing, the ones belonging to the upper part of the first branch are monotone decreasing, in both cases they all intersect once the constant solution $u\equiv 1$. Solutions of the lower part of the second branch present exactly one interior maximum point, solutions of the upper part of the second branch have exactly one interior minimum point, in both cases they have two intersections with $u\equiv 1$, and so on.

In \cite{BNW}, it was conjectured  that a similar behavior should hold also for a general nonlinearity $g$, when $p=2$. For $g$ asymptotically linear (and hence Sobolev-subcritical), this conjecture was proved to be true in \cite{ma2016bonheure}. In \cite[Corollary 1.3]{ABF} (see Corollary \ref{cor:proto}-(i) above), we prove the conjecture, without assuming any growth conditions at infinity on $g$, via shooting method.

Concerning case $p>2$, from Theorem \ref{thm:intro} we know that a non-constant solution of \eqref{Pg} arises as soon as the exponent $q>p$. Even more, Corollary~\ref{cor:proto}-(ii) guarantees that when $g(s)=s^{q-1}$, \eqref{Pg} has infinitely many solutions as soon as $q>p$. Here the eigenvalues of the operator are not involved. 
In Figure~\ref{fig:ps2}, we present some numerical simulations for the case  $p=2.1>2$. A bifurcation phenomenon from the constant solution seems to persists also when $p>2$. In this figure only the two branches of monotone solutions are detected, we refer to \cite[Section 3]{ABF} for more simulations for $p>2$. In \cite{ABF}, we conjecture that in this case infinite branches bifurcate from the same point $q=p$, giving rise to a very degenerate situation. This would be coherent with the result of Corollary~\ref{cor:proto}-(ii). We further remark that the solution found in Theorem \ref{thm:intro} is non-decreasing, so with reference to Figure \ref{fig:ps2}, it belongs to the lower (blue) branch of solutions.

\section{Proof of  Theorem \ref{thm:ABF}-(i)}\label{Sec4}
We sketch below the proof of part (i) of Theorem \ref{thm:ABF}, we refer to \cite{ABF} for more details. 

\underline{\it Step 1.} {\it (Equivalent 1-dimensional problem)}
Since we are dealing with radial positive solutions of \eqref{Pg}, we can extend $g$ to the whole of $\mathbb R$ in such a way that
$$f(s):=
\begin{cases}
g(s)-s^{p-1}\quad&\mbox{if }s\ge0,\\
0&\mbox{if }s<0
\end{cases}
$$
and write the problem in radial coordinates
\begin{equation}\label{Pr}
\begin{cases}
-(r^{N-1}|u'|^{p-2}u')'=r^{N-1}f(u)\quad\mbox{in }(0,R)\\
u'(0)=u'(R)=0.
\end{cases}
\end{equation}
We observe that while the condition $u'(R)=0$ comes from Neumann boundary conditions in \eqref{Pg}, $u'(0)=0$ is implied by symmetry and regularity of the solution. 

Then we prove (cf. \cite[Lemma 2.1]{ABF}) the following maximum principle-type result. 

\centerline{\it If $u$ solves \eqref{Pr}, then either $u>0$ in $[0,R]$ or $u\equiv -C$ for some $C\ge0$.}

As a consequence, in order to get (positive) solutions of the original problem \eqref{Pg}, it is enough to find non-constant solutions of \eqref{Pr}.

\underline{\it Step 2.} {\it (Shooting method)}
Let $\varphi_p(s):=|s|^{p-2}s$ and  $v:=r^{N-1}\varphi_p(u')$. 

We consider the ODE system
\begin{equation}\label{sys}
\begin{cases}
u'=\varphi_p^{-1}\left(\frac{v}{r^{N-1}}\right)\quad&\mbox{in }(0,R),\\
v'=-r^{N-1}f(u)&\mbox{in }(0,R),\\
u(0)=d\in[0,1],\\
v(0)=0.
\end{cases}
\end{equation}
We prove in \cite[Lemma 2.2]{ABF} global existence, uniqueness and continuous dependence for \eqref{sys}.
These results are not trivial because the system \eqref{sys} is not regular for three different reasons: at $r = 0$ we have a singularity of order $r^{-\frac{N-1}{p-1}}$ which is not integrable
when $N \ge p$;  $\varphi_p^{-1}$ is not Lipschitz continuous at 0 when $p>2$; $f$ is not Lipschitz continuous at 0 when $1<p<2$. Nevertheless, using \cite[Theorem 4]{RW97}, we are able to prove the following:

{\it \begin{itemize}
\item For all $d\in[0,1]$ there exists a unique $(u_d,v_d)$ global solution of \eqref{sys}.
\item If $d_n\to d$ then $(u_{d_n},v_{d_n})\to (u_d,v_d)$ uniformly in $[0,R]$.
\end{itemize}}
 
We observe that if $(u,v)$ solves \eqref{sys}, then $u'(0)=0$. This follows from the initial condition $v(0)=0$, cf. \cite{FLS}. Furthermore, by the definition of $v$, if $v(R)=0$, also $u'(R)=0$. Finally, for $d=0$ and $d=1$ we get the constant solutions $u\equiv 0$ and $u\equiv 1$, respectively. Hence, in order to get a non-constant solution of \eqref{Pg}, 

\centerline{\it we look for $d\in(0,1)$ such that the solution $(u_d,v_d)$ of \eqref{sys} satisfies $v_d(R)=0$.}

This procedure is referred to as {\it shooting method}.

\underline{\it Step 3.} {\it (Equivalent system in $p$-polar coordinates)}
\smallskip

\begin{minipage}[l]{.5\textwidth}
If $v(0)=v(R)=0$, by the regularity of $v$, there exists $\bar r\in(0,R)$ such that $v'(\bar r)=0$. Thus, from the equation $v'(r)=-r^{N-1}f(u)$ and by $(g_{\mathrm{eq}})$, $u(\bar r)=1$. 
\end{minipage}\;\;\;
\begin{minipage}[r]{.5\textwidth}
\includegraphics[width=.85\textwidth]{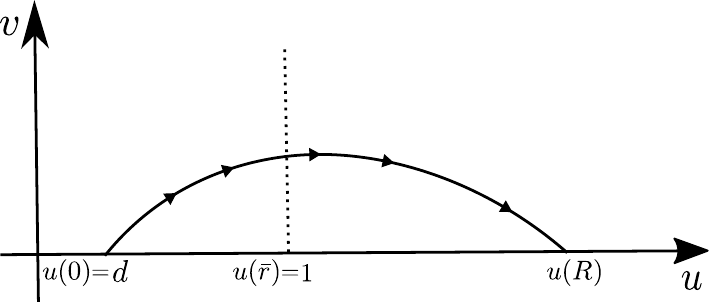}
\end{minipage}
\smallskip

Furthermore, by uniqueness, if $d\neq1$, $(u_d(r),v_d(r))\neq(1,0)$ for all $r\in[0,R]$. This means that {\it non-constant solutions of \eqref{sys} having $v(R)=0$ turn around the point $(1,0)$ in the phase plane $(u,v)$}.
 
Hence, we can pass to $p$-polar coordinates\footnote{See \cite[Section 2 and Lemma 2.3]{ABF} for the definition and properties of the functions $p$-cosine $\cos_p$ and $p$-sine $\sin_p$. Their name is due to the fact that these functions share many properties with the classical cosine and sine. For instance they are $2\pi_p$-periodic, where $\pi_p$ is the number $\pi_p=\frac{2\pi(p-1)^{1/p}}{p\sin(\pi/p)}$. Furthermore, for $p=2$, it holds $\cos_2=\cos$, $\sin_2=\sin$, and $\pi_2=\pi$. The use of these functions is common in $p$-Laplacian problems, it allows to get the equation in $\vartheta$ of the associated eigenvalue system \eqref{sys_EP} {\it not coupled} with the equation in $\varrho$.} about $(1,0)$
$$
\begin{cases}
u-1=\rho^{\frac{2}{p}}\cos_p\theta\\
v=-\rho^{\frac{2}{p'}}\sin_p\theta\end{cases}\quad\Rightarrow\quad \mbox{if } \rho>0: \quad\begin{array}{lll}u=1&\Leftrightarrow\quad \theta=(j+\frac12)\pi_p&\;(j\in\mathbb Z)\\
v=0&\Leftrightarrow\quad \theta=j\pi_p&\;(j\in\mathbb Z)
\end{array}
$$
to get the system 
\begin{equation}\label{sys_p}
\begin{cases}
\rho'(r)=\dfrac{p}{2\rho}u'\left[\varphi_p(u-1)-r^{(N-1)p'}f(u)\right]\\
\theta'(r)=r^{N-1}\left[\dfrac{p-1}{r^{(N-1)p'}}|\sin_p\theta|^{p'}+\dfrac{1}{\rho^2}(u-1)f(u)\right]\\
\theta(0)=\pi_p,\quad\rho(0)=(1-d)^{p/2}.
\end{cases}
\end{equation}
Thus, our goal becomes:

\centerline{\it Find $d\in(0,1)$ such that $\theta_d(R)=j\pi_p$ for some $j\in\mathbb Z$.}

We observe in passing that, by the equation for $\theta'$ in \eqref{sys_p} and by $(g_{\mathrm{eq}})$, we know that $\theta$ is monotone increasing. 

\underline{\it Step 4.} {\it (Using the hypothesis $0<C_1<\lambda_{k+1}^{\mathrm{rad}}$)}
By $(g_1)$ and by continuous dependence on $d$, we get for $d$ close to 1  
$$(u_d-1)f(u_d)>\left(C_1-\varepsilon\right)|u_d-1|^p=\left(C_1-\varepsilon\right)\rho_d^2|\cos_p\theta_d|^p.$$
Hence, by \eqref{sys_p}, since $C_1>\lambda_{k+1}^{\mathrm{rad}}$, for $\varepsilon>0$ sufficiently small and $d$ close to 1
$$
\begin{aligned}
\theta_d'(r)&=r^{N-1}\left[\frac{p-1}{r^{(N-1)p'}}|\sin_p\theta_d|^{p'}+\frac{1}{\rho_d^2}(u_d-1)f(u_d)\right]\\
&>r^{N-1}\left[\frac{p-1}{r^{(N-1)p'}}|\sin_p\theta_d|^{p'}+\left(C_1-\varepsilon\right)|\cos_p\theta_d|^p\right]\\
&>r^{N-1}\left[\frac{p-1}{r^{(N-1)p'}}|\sin_p\theta_d|^{p'}+\lambda_{k+1}^{\mathrm{rad}}|\cos_p\theta_d|^p\right].
\end{aligned}
$$

\underline{\it Step 5.} {\it (The associated eigenvalue problem)}
We will estimate the number of times that the solutions of the problem turn around $(1,0)$ by the number of times that the radial eigenfunctions of the Neumann $p$-Laplacian turn around $(0,0)$ in the phase plane. To this aim, we introduce the associated eigenvalue problem
\begin{equation}\label{eigenprob}
\begin{cases}
-\Delta_p\phi=\lambda^{\mathrm{rad}}|\phi|^{p-2}\phi\quad&\mbox{in }B_R,\\
\partial_\nu \phi=0&\mbox{on }\partial B_R.
\end{cases}
\end{equation}
In \cite[Theorem 1]{RW99} it has been proved what follows. 

{\it The eigenvalue problem \eqref{eigenprob} has a countable number of eigenvalues $0=\lambda_1^{\mathrm{rad}}<\lambda_2^{\mathrm{rad}}<\dots$ which go to infinity as $k\to\infty$. Furthermore, the $k$-th eigenfunction $\phi_k$ has exactly $k-1$ zeros in $(0,R)$.}

Since we are interested only in radial eigenvalues, we can write \eqref{eigenprob} as
\begin{equation}\label{EP}
\begin{cases}
-(r^{N-1}\varphi_p(\phi'))'=\lambda r^{N-1}\varphi_p(\phi)\quad\mbox{in }(0,R),\\
\phi'(0)=\phi'(R)=0.
\end{cases}
\end{equation}

We now pass to $p$-polar coordinates around $(0,0)$, that is to say
$$
\begin{cases}
\phi=\varrho^{\frac{2}{p}}\cos_p\vartheta\\
\psi:=r^{N-1}|\phi'|^{p-2}\phi'=-\varrho^{\frac{2}{p'}}\sin_p(\vartheta)
\end{cases}\Rightarrow\; \mbox{if } \varrho>0: \;\begin{array}{lll}\phi=0&\Leftrightarrow\; \vartheta=(j+\frac12)\pi_p&\;(j\in\mathbb Z)\\
\psi=0&\Leftrightarrow\; \vartheta=j\pi_p&\;(j\in\mathbb Z).
\end{array}
$$
Hence, system \eqref{EP} becomes
\begin{equation}\label{sys_EP}
\begin{cases}
\varrho'(r)=\dfrac{p}{2\varrho}\left(1-\lambda r^{(N-1)p'}\right)\varphi_p(\phi)\phi'\\
\vartheta'(r)=r^{N-1}\left[\dfrac{p-1}{r^{(N-1)p'}}|\sin_p\vartheta|^{p'}+\lambda|\cos_p\vartheta|^p\right]\\
\vartheta(0)=\pi_p,\quad\vartheta(R)=j\pi_p\quad (\exists \;j\in\mathbb Z).
\end{cases}
\end{equation}
From the second equation of \eqref{sys_EP}, we get $\vartheta'(r)>0$. 
Therefore, the fact that $\phi_{k+1}$ has exactly $k$ zeros reads as $\vartheta_{\lambda_{k+1}}(R)=(k+1)\pi_p$.

\underline{\it Step 6.} {\it (Comparing solutions with eigenfunctions)}
We now know that 
\begin{itemize}
\item[(a)] $\theta_d'(r)>r^{N-1}\left[\frac{p-1}{r^{(N-1)p'}}|\sin_p\theta_d|^{p'}+\lambda_{k+1}^{\mathrm{rad}}|\cos_p\theta_d|^p\right]$ for $d$ close to 1, by Step 4.;
\item[(b)] $\theta_d(0)=\vartheta_{\lambda_{k+1}}(0)=\pi_p$;
\item[(c)] $\vartheta_{\lambda_{k+1}}'(r)=r^{N-1}\left[\frac{p-1}{r^{(N-1)p'}}|\sin_p\vartheta_{\lambda_{k+1}}|^{p'}+\lambda_{k+1}^{\mathrm{rad}}|\cos_p\vartheta_{\lambda_{k+1}}|^p\right]$,
by Step 5.
\end{itemize}
Therefore, by Comparison Theorem 
$$\theta_d(R)>\vartheta_{\lambda_{k+1}}(R)=(k+1)\pi_p\quad \mbox{as }d\sim1,$$
that is to say, the solution performs more than $k$ half-turns around $(1,0)$ in the phase plane.
Then, by 
continuous dependence of $(u_d,v_d)$ and hence of $(\rho_d,\theta_d)$ on $d$, and by the fact that $\theta_0(R)=\pi_p$ (i.e., 0 turns), we obtain that there exist $d_1,\dots,d_k\in(0,1)$  such that

\begin{minipage}[l]{.45\textwidth}
$\theta_{d_j}(R)=(j+1)\pi_p$ for any $j=1,\dots,k$,
\end{minipage}\;\;\;
\begin{minipage}[r]{\textwidth}
\smallskip
\includegraphics[width=.5\textwidth]{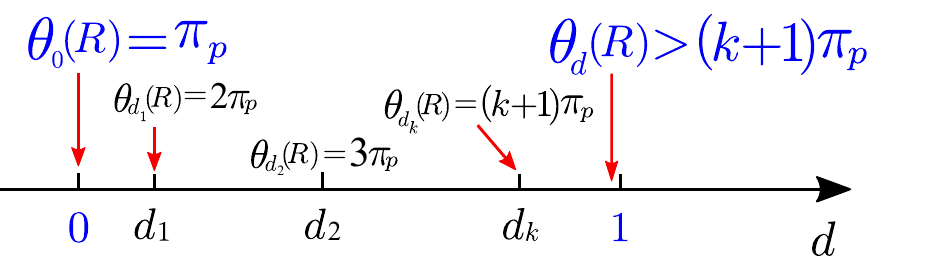}\smallskip
\end{minipage}
which correspond to the $k$ non-constant radial solutions $u_1,\dots,u_k$. Furthermore, since $\theta_{d_j}(0)=\pi_p$, $\theta_{d_j}(R)=(j+1)\pi_p$, and $\theta_{d_j}$ is monotone increasing, we immediately get that $u_j-1$ has exactly $j$ zeros for any $j=1,\dots,k$. 

\begin{remark}
With reference to Figure \ref{fig:p2}, we observe that in the pure power case $g(s)=s^{q-1}$, from Step 2. we can see that the solutions detected in Theorem \ref{thm:ABF} belong to the lower parts of the branches, since they all satisfy $u(0)=d<1$. 
\end{remark}

\section{The case $1<p<2$}\label{Sec5}
Al already mentioned in the Introduction, the case $C_1=0$ corresponds to the case $1<p<2$ for the prototype nonlinearity $g(s)=s^{q-1}$, $q>p$. This is the reason why in this section, devoted to the case $1<p<2$, we start with some comments on the proof of Theorem \ref{thm:ABF}-(iii).
The proof of this part is rather technical, we want to highlight here only the main differences with the case $C_1\neq 0$ which are responsible for the surprising result found. 
To this aim, we start observing that in the proof of part (i) it is crucial to have an estimate of the number of times that the solution of \eqref{sys}, shot from a point $d$ close enough to 1 of the $u$-axis, turns around the point $(1,0)$ in the phase plane: in Step 6. of the previous section we end up with the following estimate from below $\theta_d(R)>(k+1)\pi_p$ for $d\sim1$. Instead, in the case $C_1=0$, thanks to an adaptation of \cite[Corollary 5.1]{BZ} (see \cite[Lemma 2.8]{ABF}), we get the following result:  

\centerline{\it If $R > R_*(k)$, there exists $d_*\in(0,1)$ such that $\theta_{d_*}(R)>(k+1)\pi_p$.}

This means that we know that the number of half-turns is greater than $k+1$ for solutions shot at the {\it finite} distance $1-d^*$ from the point $(1,0)$, and not in the limit as $d\to1$. Furthermore, by $(g_1)$ and Gronwall's inequality, we prove for $\lambda=C_1=0$ 
$$\theta_d(R)\to\vartheta_0(R)=\pi_p\quad\mbox{\it as }d\to1.$$ 
This allow us to make the continuous-dependence procedure effective both for solutions shot from $u(0)=d\in (0,d^*)$ and for solutions shot from $u(0)=d\in (d^*,1)$. In this way, we obtain the double of the solutions found for $C_1\in(0,\infty)$, as represented in the following picture. \smallskip

\begin{center}
\includegraphics[width=.65\textwidth]{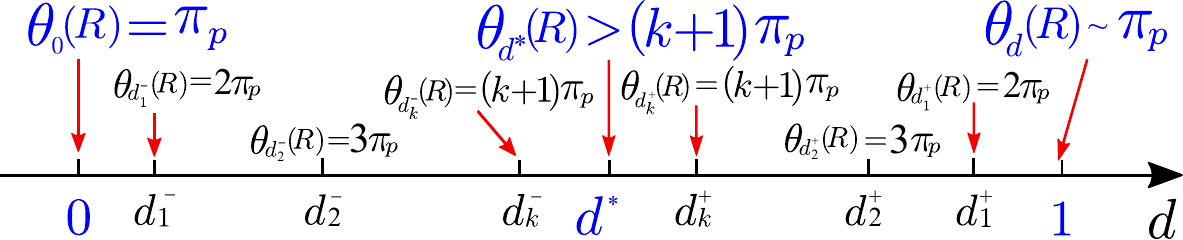}
\end{center}
More precisely, from one side, by continuous dependence on $d$ and since $\theta_0(R)=\pi_p$, we have that
$$\mbox{\it there exist }\; d^-_1, \dots,d_k^-\in(0,d_*)\;\;\mbox{\it s.t. }\theta_{d^-_j}(R)=(j+1)\pi_p\;\;\mbox{\it  for all }j=1,\dots,k.$$
On the other side, again by continuous dependence on  $d$, being $\theta_{d_*}(R)>(k+1)\pi_p$, we obtain
$$\mbox{\it there exist }\;d^+_1, \dots,d_k^+\in(d_*,1)\;\;\mbox{\it s.t. }\theta_{d^+_j}(R)=(j+1)\pi_p\;\; \mbox{\it  for all }j=1,\dots,k.$$
\smallskip

Some comments are now in order. 
In \cite[Section 3]{ABF} some numerical simulations performed for $N=1$, $R=1$, $p=1.97<2$, and $g(s)=s^{q-1}$ show that for values of $p<2$ sufficiently close to $2$ the branches of solutions persist. 
\begin{figure}[h!t]
\begin{center}
\includegraphics[width=8cm]{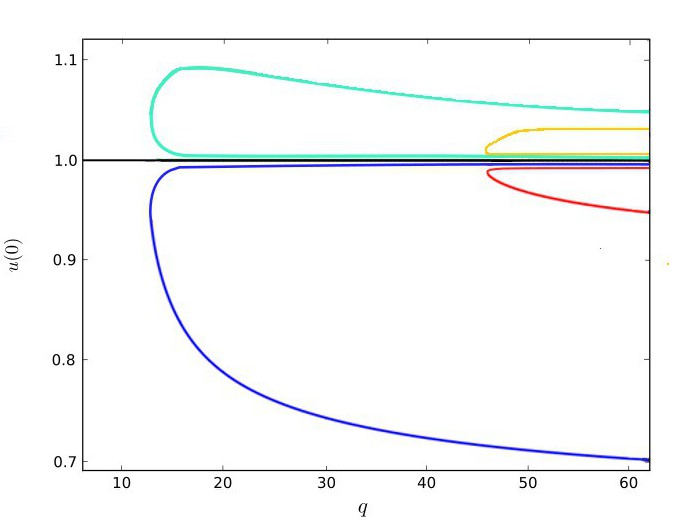}
\caption{Qualitative representation of the first four branches of non-constant solutions for problem \eqref{Pg} in the case $N=1$, $1 \ll p<2$, $R=1$, $g(s)=s^{q-1}$.}\label{fig:pi2}
\end{center}
\end{figure}
\noindent Differently from what we found for $p=2$, now each branch splits into two and both the upper and the lower part of the branches {\it fold}, as represented in Figure \ref{fig:pi2}. This heuristically explains why for $p<2$ we find the double of solutions with respect to the case $p=2$. Indeed, the shape of the branches is coherent with the result found in Corollary \ref{cor:proto}-(iii), since for every value of $q>p$, each {\it folded} branch contains now {\it two} different solutions having the same oscillatory behavior. Furthermore, none of the branches seem to bifurcate from the constant solution $u\equiv1$, but each of them seem to converge to the constant solution as $q\to\infty$. It looks like as the bifurcation point has escaped to infinity. 

\section*{Acknowledgments}
The authors wish to thank Prof. Alberto Parmeggiani for interesting suggestions on future perspectives about elliptic domains, and Dr. Alberto Boscaggin for raising a useful question about the pre-existing condition $(g_0)$ on page 3, which led us to generalize it into the present form.
B. Noris acknowledges the support of the project ERC Advanced Grant 2013 n. 339958: ``Complex Patterns for Strongly Interacting Dynamical Systems --
COMPAT''. F. Colasuonno was supported by the INdAM - GNAMPA Project 2017 ``Regolarit\`a delle soluzioni viscose per equazioni a derivate parziali non lineari degeneri" and by University of Bologna, funds for selected research topics. 

\bibliographystyle{alpha}

\end{document}